\documentclass[12pt]{amsart}
\usepackage{amsthm, amsxtra, latexsym, mathrsfs, amscd}
\usepackage{graphicx}
\usepackage{epsfig}
\usepackage{amsfonts}
\usepackage{amscd}
\usepackage{amsfonts}
\usepackage{psfrag}
\usepackage{pstricks}
\usepackage{mathdots}
\usepackage{psfrag}
\usepackage{amsthm,amsmath,amssymb,graphics,hyperref}
\usepackage{amsmath,amssymb,hyperref}
\usepackage{enumerate}
\usepackage[titletoc,title]{appendix}
\input xy
\xyoption{all}

\theoremstyle{plain}
\newtheorem{theorem}{Theorem}[section]
\newtheorem{proposition}[theorem]{Proposition}
\newtheorem{lemma}[theorem]{Lemma}
\newtheorem{corollary}[theorem]{Corollary}
\theoremstyle{definition}

\newcommand{\R}{\mathbb{R}}
\newcommand{\bc}{\mathbb{C}}
\newcommand{\ra}{{\rightarrow}}

\let\cal\mathcal

\def\P{\mathbb{P}}

\begin{document}
\title[Entropy]
        { Topological Entropy and bulging deformation of real projective structures on surface }
        \author{Patrick Foulon and Inkang Kim}
        \date{}
        \maketitle

\begin{abstract}
In this paper we study  the deformation of strictly convex real
projective structures on a closed surface. Specially we study the deformation in terms of the entropy
 on bulging deformations.
As a byproduct we construct a sequence of divergent structures whose topological
entropy converges to a designated number between 0 and 1.
\end{abstract}
 \footnotetext[3]{I. Kim
gratefully acknowledges the partial support of  Grant
(NRF-2014R1A2A2A01005574)  and a warm support of CIRM during his
visit.}
\section{Introduction}
A real projective structure on a closed surface $S$ is a maximal atlas where the transition functions are the restriction of elements in
$PGL(3,\R)$. The developing map of the real projective structure defines a local homeomorphism from $\widetilde S$ into $\R\P^2$.
When the image of the developing map is a convex domain, the real projective structure is said to be convex.
In this note, we are interested in a strictly convex real projective structure. A strictly convex domain in $\R\P^2$ can be considered as a deformation of a conic where the structure is real hyperbolic. In this sense, the moduli space of real projective structures contain the Teichm\"uller space. Indeed it is known that the moduli space is a holomorphic vector bundle over the Teichm\"uller space with the fibres being holomorphic cubic differentials \cite{La,Lo2}. Furthermore it is a K\"ahler manifold with the Teichm\"uller space being
totally geodesic equipped with Weil-Petersson metric \cite{KZ}.

In the language of Higgs bundles, Hitchin \cite{Hi} singled out a component in a character variety of $SL(3,\R)$. This Hitchin component is the set of holonomy representations of strictly convex real projective structures. Labourie \cite{La1} also characterized this component as Anosov representations. But until now, its compactification and dynamical properties remain mysterious.

In this note, we analyze the dynamical properties of geodesic flow in terms of the entropy with respect to some invariant measures. By Crampon \cite{crampon} it is known that the topological entropy of the geodesic flow is less than or equal to  that of hyperbolic structure, and the equality occurs only when the real projective structure is a hyperbolic structure.
Later Nie \cite{Nie} constructed a sequence of real projective structures whose topological entropy tends to zero.
More recently, Zhang \cite{Zhang} constructed a sequence of divergent real projective structures whose topological entropy tends to a prescribed number between 0 and 1. He uses so-called internal parameters to construct such a sequence. In this paper, we use different parameter, bulging parameter or vertical parameter, to construct such a sequence.
 Goldman \cite{Goldman} gave parameters for strictly convex real projective structures as follows. Take a pants decomposition of $S$. For each pair of pants, there are boundary parameters, two parameters for each boundary curve. There are two more internal parameters to determine a projective structure on a pair of pants. Hence 8 parameters for each pair of pants. To glue two pairs of pants, there are two parameters, called a horizontal twisting parameter, and a vertical parameter (or bulging parameter). In this paper, we are concerned with a vertical parameter.

We study a new deformation, originally due to Goldman \cite{Go}, called a bulging deformation.
If $\rho:\pi_1(S)\ra SL(3,\mathbb R)$ is a holonomy representation of a strictly  convex real projective structure on a closed surface $S=\Omega/\rho(\pi_1(S))$ of genus at least 2, then for every closed loop $\gamma$ on $S$, considered as a conjugacy  class in $\pi_1(S)$, $\rho(\gamma)$ can be diagonalized with eigenvalues
$\lambda_1>\lambda_2>\lambda_3$  so that the corresponding  eigenvectors
represent  attracting,  neutral, repelling fixed points $\gamma_+,\gamma_0,\gamma_-$ of $\rho(\gamma)$ on $\mathbb{RP}^2$. For the definition of bulging deformation along $\gamma$, see Section \ref{bulging}. Intuitively, bulging along $\gamma$ is stretching the domain $\Omega$ to the direction of $\gamma_0$, while  the earthquake along $\gamma$ is sliding the domain along the axis $\overline{\gamma_-,\gamma_+}$ of $\rho(\gamma)$.

Our main theorem can be summarized as:
\begin{theorem}\label{cylinder}(1) Suppose $\gamma$ is separating. If we let the vertical  parameter $s$ go to $\infty$, the domain converges, in Gromov-Hausdorff topology with a base point,  to a properly convex domain (but not
strictly convex), and the strictly convex structure on $S$,  with a base point  on the left component $S'$ of
 $S\setminus \gamma$, converges
 to a convex structure (but not strictly convex)  of infinite volume.  In this case, the projective structure converges
 to a projective structure on $S'$ with  almost Euclidean half open cylinder attached
 along $\gamma$. The triangle determined by
 $\gamma_-,\gamma_+,\gamma_0$ projects onto this half open cylinder.

 If $s\ra -\infty$, then the domain converges to the one with the
 axis $\overline{\gamma_-,\gamma_+}$ of $\gamma$ on its boundary. In this case, the projective
 structure converges to a projective structure on $S'$ with a cylinder corresponding to
 $\gamma$. A neighborhood of $\overline{\gamma_-,\gamma_+}$ projects
 to the cylinder.

A similar statement holds for the right component of $S\setminus\gamma$ if we choose a base point on the right component.

(2) If $\gamma$ is non-separating, the limit structure is the one obtained by attaching cylinders along each copy of $\gamma$.
\end{theorem}

{One can do bulging along a disjoint union of simple closed curves consecutively. Such a maximal union of disjoint, non-parallel simple closed curves is a pants decomposition of $S$.}
As an application we get
\begin{theorem}Let $\rho_s:\pi_1(S)\ra SL(3,\R)$ be a smooth family of bulging deformation.   If $\rho_s$ represents a family of bulging deformations on pants decomposition $P$, then the topological entropy tends to the maximum of topological entropies on pairs of pants $S\setminus P$ as $s\ra\infty$. Furthermore one can construct a divergent sequence whose topological entropy tends to a prescribed number between 0 and 1.
\end{theorem}
One can consider other entropies using different invariant measures. A natural geodesic flow invariant measure other than Bowen-Margulis measure is Sinai-Ruelle-Bowen measure. In another paper by the authors we deal with this measure entropy using SRB measure \cite{FK}.

\section{Preliminaries}
\subsection{Projective structure}
Let $\Omega$ be a strictly convex domain with $C^1$ boundary equipped with a Hilbert metric. The corresponding Finsler norm on $T\Omega$ is denoted by $F$. If it admits a compact quotient manifold $M=\Omega/\Gamma$, $\Gamma$ is Gromov hyperbolic and the geometry behaves like a negatively curved case.
More generally, if a convex domain $\Omega$ is not a triangle and admits a finite volume quotient, then it is strictly convex \cite{Marquis}.

For a given $w=(x,[\xi])\in H\Omega=T\Omega\setminus \{0\}/\R_+^*$, the unstable manifold $W^u$ passing through $w$ is defined to be
$$W^u(w)=\{(y,[\phi])\in H\Omega| \xi(-\infty)=\phi(-\infty), y\in \cal H_{\sigma w}\}.$$ Here $\xi(-\infty)$ denotes $\gamma_\xi(-\infty)$ where $\gamma_\xi$ is the geodesic determined by $\xi$, and $\sigma w=(x,[-\xi])$ is a flip map, $\cal H_w$ is the horosphere based at $\xi(\infty)$ passing through $x$.
Similarly one can define a stable manifold
$$W^s(w)=\{(y,[\phi])\in H\Omega| \xi(\infty)=\phi(\infty), y\in\cal H_w\}.$$ These stable and unstable manifolds are $C^1$ if $\partial \Omega$ is $C^1$.

The tangent spaces of $W^u$ and $W^s$ form unstable and stable vectors in $TH\Omega$, i.e., along the geodesic flow, they
 expand or decay exponentially. It is known \cite{Be} that the geodesic flow on $HM$ is Anosov with invariant decomposition
$$THM=\R X\oplus E^s\oplus E^u,$$ where $X$ is the vector field generating the geodesic flow.
\subsection{Hilbert metric}\label{Hilbertmetric}
Suppose that $\Omega$ is a (not necessarily strictly) convex domain in
$\mathbb{RP}^n$. For $x\neq y\in \Omega$, let $p,q$ be the
intersection points of the line $\overline{xy}$ with $\partial\Omega$ such that $p,x,y,q$ are in this order. The
Hilbert distance is defined by
$$d_\Omega(x,y)=\frac{1}{2}\log \frac{|p-y||q-x|}{|p-x||q-y|}$$ where $| \cdot |$ is a Euclidean
norm. This metric coincides with the hyperbolic metric if
$\partial\Omega$ is a conic. The Hilbert metric is Finsler rather than
Riemannian. The Finsler norm $F=||\cdot||$ is given, for $x\in\Omega$ and a
vector $v$ at $x$, by
$$||v||_x=\frac{1}{2}\big(\frac{1}{|x-p^-|}+\frac{1}{|x-p^+|}\big)|v|$$ where $p^\pm$
are the intersection points of the line with $\partial\Omega$, defined by $x$ and $v$ with the
obvious orientation, and where again $| \cdot |$ is
a Euclidean norm. Then it is classical that the metric induced by
this Finsler norm is the Hilbert metric. In particular, it is reversible.

Choose an affine set $A$ containing $\Omega$ with a Euclidean norm
$| \cdot |$. Let $\mathrm {Vol}$ be a Lebesgue measure on $A$ normalized by
$\mathrm{Vol}_x(\{v\in A:|v|<1\})=1$. { Here $T_xA$ is naturally identified with $A$.} Then for any Borel set $\cal A\subset
\Omega\subset A$, one can define a measure
\begin{eqnarray}\label{Finslervolume}
\mu_\Omega(\cal A)=\int_{\cal A} \frac{d\mathrm {Vol}_x}{\mathrm {Vol}(B_x(1))}
\end{eqnarray}
where $B_x(1)=\{v\in T_x\Omega: ||v||_x<1\}$. This measure turns
out to be the Hausdorff measure induced by the Hilbert metric
\cite{BBI}.

From the definition it is clear that for two convex domains
$\Omega_1\subset \Omega_2$,
\begin{enumerate}
\item  $||v||_x^{\Omega_2}\leq ||v||_x^{\Omega_1}$;
\item $d_{\Omega_2}(x,y)\leq d_{\Omega_1}(x,y)$;
\item  $B_x^{\Omega_1}(1)\subset B_x^{\Omega_2}(1)$;
\item  $\mu_{\Omega_2}(\cal A)\leq \mu_{\Omega_1}(\cal A)$ for any
Borel set $\cal A$.
\end{enumerate}

\subsection{Invariant measures}
For a geodesic flow on $H\Omega$, the maximal entropy $h_\mu(\phi)$ of the probability measure $\mu$ is known to be realized at the Bowen-Margulis measure.
This entropy is equal to the topological entropy $h_{top}(\phi)$  of the geodesic flow $\phi$ and also it is equal to the exponential growth of the lengths of  closed geodesics:
$$\lim_{R\ra\infty} \frac{\log \#\{[\gamma]|\ell(\gamma)\leq R\}}{R}.$$  This is again equal to the critical exponent of the associated Poincar\'e series.

\section{Bulging deformation}\label{bulging}
\subsection{Definition of bulging}
For a closed strictly convex projective surface, every element in
the holonomy group is hyperbolic, i.e., conjugate to a diagonal
matrix $D(\alpha_1,\alpha_2,\alpha_3)$ with mutually distinct
$\alpha_i$.

 Let $S$ be a closed
 surface of genus at least 2 and $\gamma$  a closed loop. Then $\pi_1(S)=\Gamma_1*_{\langle \gamma \rangle} \Gamma_2$ or HNN extension depending on whether $\gamma$ is separaing or not.     After
conjugation, we can assume that
$$\rho(\gamma)=\gamma_{t_0}=\left [\begin{matrix}
  e^{t^0_1} & 0 & 0 \\
  0   & e^{t^0_2} & 0 \\
  0   & 0 & e^{t^0_3}
  \end{matrix}\right],\ t^0_1>t^0_2>t^0_3,\ t^0_1+t^0_2+t^0_3=0$$ in $SL(3,\R)$ where $\rho:\pi_1(S)\ra SL(3,\R)$ is the holonomy representation of a strictly convex real projective structure on the surface. Let ${\gamma_t=\left[\begin{pmatrix}
e^{tt_1^0} & 0 & 0 \\
  0   & e^{tt_2^0} & 0 \\
  0   & 0 & e^{tt_3^0}
  \end{pmatrix}\right],\ tt_1^0>tt_2^0>tt_3^0, t_1^0+t_2^0+t_3^0=0}
       $ be the one-parameter group generated by
  $\rho(\gamma)=\gamma_{t_0}$. If $\Omega$ is a strictly convex domain  so
that $S=\Omega/\rho(\pi_1(S))$, the eigenspaces $\R v_1,\R v_2,\R
v_3$ corresponding to $e^{t^0_1},e^{t^0_2},e^{t^0_3}$ respectively, define three
points $\gamma_+, \gamma_0, \gamma_-$ in $\mathbb{RP}^2$ where $\gamma_\pm$
are attracting, repelling fixed points of $\rho(\gamma)$ on
$\partial \Omega$, and $\gamma_0$ is outside  $\Omega$. Consider a
triangle $\triangle$ with  left
and right vertices corresponding to $\gamma_+,\gamma_-$, top vertex
to $\gamma_0$. Then the dynamics of $\rho(\gamma)$ on $\triangle$ is
from the right vertex to the left and top vertices, from the top
vertex to the left vertex. Inside $\triangle$ the orbits of
$\gamma_t$ are arcs tangent to $\triangle$ Let  $C$ be the arc of $\partial\Omega$ from $\gamma_-$ to
$\gamma_+$.  We naturally orient $\gamma$ by the orientation from $\gamma_-$ to $\gamma_+$.

The earthquake map of $S$ along $\gamma$ is given by the right Dehn
twist along $\gamma$, which amounts to moving right hand side of the
lifts of $\gamma$ by the amount $t$ in $\Omega$. This can be
realized by conjugating the action of the right side of $S\setminus\gamma$ (if
$\gamma$ is separating) by
$$\tau_t=\left [\begin{matrix}
  e^{t} & 0 & 0 \\
  0   & 1 & 0 \\
  0   & 0 & e^{-t}
  \end{matrix}\right],$$ and correspondingly using HNN
extension for non-separating case.

Obviously this earthquake deformation does not change the domain
$\Omega$ if the starting domain is conic. To deform the domain, we do so-called {\it bulging deformation} originated from  \cite{Go}.
We want to replace $C$ by another curve tangent to $\triangle$. This can be realized by
conjugating the right side of $S\setminus\gamma$ (if $\gamma$ is
seperating) by $$O_s=\left[\begin{matrix}
                   e^{-\frac{1}{3}s} & 0 & 0 \\
                   0 & e^{\frac{2}{3}s} & 0 \\
                   0 & 0 & e^{-\frac{1}{3}s}\end{matrix}\right]. $$  {For  $\gamma$ non-separating case,
$$\rho(\pi_1(S))=\Gamma=\Gamma_1*_{\langle\gamma\rangle}=\langle \Gamma_1,
\gamma_2\rangle\,,$$
 where  $\gamma_2\in\Gamma\setminus\Gamma_1$ conjugates two subgroups of
$\Gamma_1$ that are isomorphic to $\langle\gamma\rangle$. The bulging
deformation gives rise to
$$\langle \Gamma_1,
O_s\gamma_2\rangle.$$  Note here that since $S$ is oriented, once $\gamma$ is oriented, the left and right side of $\gamma$ make sense locally.}

 What $O_s$ does to the
domain is stretching the right  side of $\Omega\setminus \tilde\gamma$ to $\gamma_0$ direction, which entails to
move the boundary arc $C$ to $O_s C$
outside $\Omega$ if $s>0$, to one inside $\Omega$ if $s<0$. Here $\tilde\gamma$ is a lift of $\gamma$ whose ends are $\gamma_-$ and $\gamma_+$ and $\triangle$ is a triangle whose vertices are $\gamma_\pm$ and $\gamma_0$. Note that under this stretching, the deformed part is contained in $\triangle$.

\begin{lemma}Bulging the right side of $\tilde\gamma$ is projectively equivalent to debulge the left side of $\tilde\gamma$.
\end{lemma}
\begin{proof}If we apply $O_s^{-1}$ to the deformed domain, the right side will move back
to the original one but the left side will be debulged, i.e., the left side will be shrunk toward
$\tilde\gamma$.
\end{proof}

 One does this operation to each lift of $\gamma$ to obtain a deformed domain $\Omega_s$. {More precisely} one can obtain the domain $\Omega_s$ as follows.

Let $\{\text{lifts of }\gamma\}=\cup \alpha_i \tilde\gamma$ where $\alpha_i\in \pi_1(S)$.
Suppose first that $\gamma$ is separating and with a natural orientation on $\gamma$ that $S\setminus\gamma$ has left and right components $S^L$, $S^R$ respectively.
Take   a component $K^L$ of $\Omega\setminus \{\text{lifts of }\gamma\}$ bounded by $\tilde\gamma$, which projects down to $S^L$ and an adjacent   component $K^R$ which projects down to $S^R$.  Note that {$\overline{K^L}\cap \overline{K^R}=\tilde\gamma$}.
Let $(K^L\subset) \Omega^L$ be the left component of $\Omega\setminus\tilde\gamma$ and similarly $(K^R\subset)\Omega^R$ be the right component.
Start doing bulging along the { lifts of $\gamma$ which are the sides}  of $K^L$.
The bulging along $\tilde\gamma$ is stretching $\Omega^R$ toward $\gamma_0$ and $\Omega^L$ unchanged. Note that the deformed domain under this first bulging along $\tilde\gamma$ is contained in $\Omega\cup\triangle$. Since $\gamma$ is separating, all the lifts of $\gamma$ which bound $K^L$ are oriented
in the way that if we walk along the curve in the positive direction, then $K^L$ lies on the left side. See Figure 1.

\begin{figure}
\begin{center}
%
%



\scalebox{0.9} 
{
\begin{pspicture}(3,-1.238125)(17.881874,2.238125)
\psellipse[linewidth=0.04,dimen=outer](7.08,0.6296875)(1.66,1.47)
\psline[linewidth=0.04cm](7.44,-0.8003125)(14, 0.6)
\psline[linewidth=0.04cm](7.44,2.07)(14, 0.6)
\psline[linewidth=0.04cm]{<-}(7.44,2.07)(7.46,-0.8003125)
\rput(7.9, 1){$K^R$}
\rput(7.7, 0){$\tilde\gamma$}
\rput(7.5,2.3){$\gamma_+$}
\rput(7.5,-1.2){$\gamma_-$}
\rput(10, 0.3){$\triangle$}
\psline[linewidth=0.04cm]{->}(8.2,1.7)(8.6,0.08)
\psline[linewidth=0.04cm]{->}(7.0,2.07)(5.5,1)
\psline[linewidth=0.04cm](5.8,1.9)(5.5,1)
\psline[linewidth=0.04cm](5.8,1.9)(7,2.07)
\psline[linewidth=0.04cm]{<-}(7.0,-0.8)(5.46,0.5)
\psline[linewidth=0.04cm](7,-0.8)(5.5,-0.6)
\psline[linewidth=0.04cm](5.46,0.5)(5.5,-0.6)
\rput(5.0,-0.6){$\alpha_i\triangle$}
\rput(6.7, -0.193){$\alpha_i(\tilde\gamma)$}
\rput(6.5, 0.8){$K^L$}
\psarc[linestyle=dashed](7.44,0.63484375){1.43515625}{-90.0}{90.0}
\usefont{T1}{ptm}{m}{n}

\rput(14.2,0.5){$\gamma_0$}
\end{pspicture}
}
\caption{Bulging deformation and the limit domain $\Omega_\infty$ when $\gamma$ is separating}
\end{center}
\end{figure}
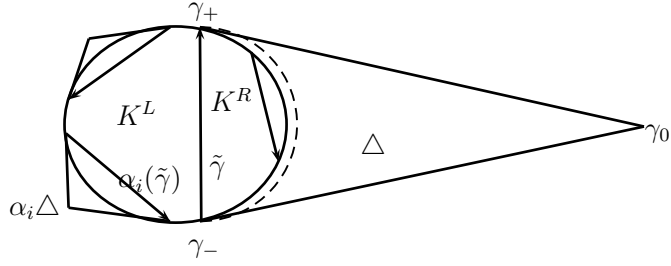

 Hence the bulging along the
side $\alpha_i\tilde\gamma$ of $K^L$ is stretching the right side of $\Omega\setminus \alpha_i\tilde\gamma$ toward $\alpha_i\gamma_0$. Then the deformed domain is contained in
$\Omega_\infty^L=\Omega \cup \cup_i \alpha_i\triangle$. Note that along the deformation,
$K^L$ remains unchanged.  Now we move onto the image $O_s K^R$ of $K^R$ to do bulging.  The boundary
curve of $K^R$ is oriented that if we walk along the curve in the positive direction, $K^R$ is on the right side. Hence when we do bulging along the side of $O_s K^R$, we stretch the domain toward $O_s K^R$. If we want
$K^L \cup O_s K^R$ unchanged, doing bulging along the side $l$ of $O_s K^R$ is projectively equivalent to  debulging
the left hand side of  $l$.  Hence the deformed domain still remains in $\Omega_\infty^L$.
In this way we can conclude that the bulged domain $\Omega_s$ staisfies
\begin{eqnarray}\label{sepL}
K^L=K_s^L\subset  K^L\cup O_s K^R \subset \Omega_s \subset \Omega_\infty^L.
\end{eqnarray}
\begin{lemma}\label{1} Suppose $\gamma$ is separating.
The Hausdorff limit  of $\Omega_s$  with a base point in $K^L$ is $\Omega^L_\infty$.
\end{lemma}
\begin{proof}Pick any point $p\in\triangle$. Then there exists $s_0$ such that for any $s>s_0$,
{$p\in O_s K^R$}.
This can be easily seen  by choosing a line from $\gamma_0$ passing through $p$ and intersecting $\tilde\gamma$. $O_s$ stretches this line toward $\gamma_0$, and if $s$ is large enough, the small segment near $\tilde\gamma$ will be stretched to a segment containing $p$. Since
bulging or debulging occurs outside $K^L\cup O_sK^R$ as explained above, $K^L\cup O_sK^R$ remains
unchanged during the bulging process, hence $p\in K^L\cup O_sK^R\subset \Omega_s$. Since $p$ is an arbitrary point of $\triangle$, $\triangle$ is contained in the Hausdorff limit of $\Omega_s$.
By the same reasoning, each triangle $\alpha_i\triangle$ is contained in the Hausdorff limit of $\Omega_s$. Hence the Hausdorff limit of $\Omega_s$ with a base point in $K^L$ is equal to $\Omega_\infty^L$.
\end{proof}
There is another way to describe the bulged domain $\Omega_s$ containing $K^R$. Instead of starting from $K^L$, one can   start bulging from $K^R$.  As already noticed, doing bulging along the side $l$ of
$K^R$ is stretching the right side of $l$ toward $K^R$. If we want $K^R$ unchanged, this process
is projectively equivalent to debulge the left side of $l$.  If we do the same reasonning as above,
we find that
\begin{eqnarray}\label{sepR}
K^R=K_s^R \subset \Omega_s \subset \Omega
\end{eqnarray}
\begin{lemma}\label{2} Suppose $\gamma$ is separating.
In this case, the Hausdorff limit of $\Omega_s$ with a base point in $K^R$ is equal to $K^R$.
\end{lemma}

Now suppose that $\gamma$ is non-separating. In this case the situation is a little bit complicated since
all the lifts of $\gamma$ which bound $K^L$ are not oriented
in the way that if we walk along the curve in the positive direction, then $K^L$ lies on the left side
as in the separating case. For some lift of $\gamma$, $K^L$ lies on the right side when we walk along the curve in the positive direction. We start bulging from $K^L$. Bulging along $\tilde\gamma$ is the same as before, hence the deformed domain is contained in  $\Omega\cup\triangle$. When we walk along the oriented boundary curve $l$ of $K^L$, if
$K^L$ is on the left side we bulge right side, hence the bulged part is contained in the attached
triangle, and if $K^L$ is on the right side, we debulge the left side, hence the debulged part is
contained in $\Omega$. See Figure 2.

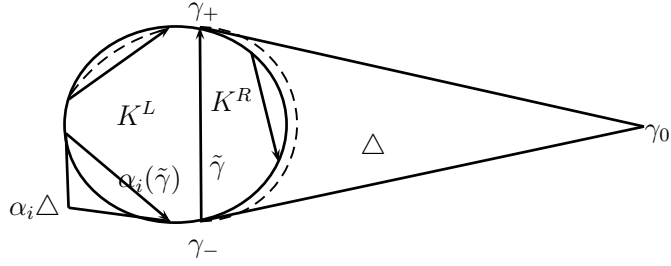
\begin{figure}
\begin{center}
%
%



\scalebox{0.9} 
{
\begin{pspicture}(3,-1.238125)(17.881874,2.238125)
\psellipse[linewidth=0.04,dimen=outer](7.08,0.6296875)(1.66,1.47)
\psline[linewidth=0.04cm](7.44,-0.8003125)(14, 0.6)
\psline[linewidth=0.04cm](7.44,2.07)(14, 0.6)
\psline[linewidth=0.04cm]{<-}(7.44,2.07)(7.46,-0.8003125)
\rput(7.9, 1){$K^R$}
\rput(7.7, 0){$\tilde\gamma$}
\rput(7.5,2.3){$\gamma_+$}
\rput(7.5,-1.2){$\gamma_-$}
\rput(10, 0.3){$\triangle$}
\psline[linewidth=0.04cm]{->}(8.2,1.7)(8.6,0.08)
\psline[linewidth=0.04cm]{<-}(7.0,2.07)(5.5,1)
\psline[linewidth=0.04cm]{<-}(7.0,-0.8)(5.46,0.5)
\psline[linewidth=0.04cm](7,-0.8)(5.5,-0.6)
\psline[linewidth=0.04cm](5.46,0.5)(5.5,-0.6)
\rput(5.0,-0.6){$\alpha_i\triangle$}
\rput(6.7, -0.193){$\alpha_i(\tilde\gamma)$}
\rput(6.5, 0.8){$K^L$}
\psarc[linestyle=dashed](7.44,0.63484375){1.43515625}{-90.0}{90.0}
\psarc[linestyle=dashed](7.5, -0.3){2.4}{105.0}{145.0}
\usefont{T1}{ptm}{m}{n}

\rput(14.2,0.5){$\gamma_0$}
\end{pspicture}
}
\caption{Bulging deformation and the limit domain $\Omega_\infty$ when $\gamma$ is non-separating}
\end{center}
\end{figure}

 In this way, we can either attach a triangle, or use $\Omega$ to obtain a limit
domain $\Omega_\infty^L$.  We do the same thing starting from $K^R$. This way we  have
\begin{eqnarray}\label{nonsep}
K^L=K_s^L\subset \Omega_s \subset \Omega_\infty^L,\
K^R=K_s^R\subset \Omega_s \subset \Omega_\infty^R.\end{eqnarray}
\begin{lemma}\label{nonseparating}Suppose $\gamma$ is non-separating.
The Hausdorff limit domain will be the one by attaching a triangle when one has to bulge, or by
cutting out the domain when one has to debulge.
\end{lemma}

\subsection{Limit structure of bulging}
Note that $\tau_t$ and $O_s$ commute and
$$\tau_tO_s=\left[\begin{matrix}
              e^{t-\frac{1}{3}s} & 0 & 0 \\
              0      &  e^{\frac{2}{3}s} & 0 \\
              0  & 0 &     e^{-t-\frac{1}{3}s}\end{matrix}\right]$$
hence, generate whole 2-dimensional abelian group $A$ in the Iwasawa
decomposition $SL(3,\R)=KAN$. In Goldman's notation \cite{Goldman}
for Fenchel-Nielsen coordinate of projective structures, $\tau_t$
corresponds to horizontal twisting parameter, and $O_s$ corresponds
to vertical twisting parameter.

 We will denote the representation
$\rho(\Gamma_1)*_{\langle \gamma \rangle} (\tau_t O_s)\rho(\Gamma_2)(\tau_t O_s)^{-1}$ by $\rho_{t,s}$ and
corresponding projective structure by $S_{t,s}=\Omega_{t,s}/\rho_{t,s}$. When $t=0$, we set $\Omega_{0,s}=\Omega_s, \rho_{0,s}=\rho_s$ for simplicity.

\begin{theorem}\label{cylinder}(1) Suppose $\gamma$ is separating. If we let the vertical  parameter $s$ go to $\infty$, the domain converges, in Gromov-Hausdorff topology with a base point,  to a properly convex domain (but not
strictly convex), and the strictly convex structure on $S$,  with a base point on on the left component $S'$ of
 $S\setminus \gamma$, converges
 to a convex structure (but not strictly convex)  of infinite volume.  In this case, the projective structure converges
 to a projective structure on $S'$ with  almost Euclidean  \marginpar {}\\half open cylinder attached
 along $\gamma$. The triangle determined by
 $\gamma_-,\gamma_+,\gamma_0$ projects onto this half open cylinder.

 If $s\ra -\infty$, then the domain converges to the one with the
 axis $\overline{\gamma_-,\gamma_+}$ of $\gamma$ on its boundary. In this case, the projective
 structure converges to a projective structure on $S'$ with a cylinder corresponding to
 $\gamma$. A neighborhood of $\overline{\gamma_-,\gamma_+}$ projects
 to the cylinder. {In terms of representations, this case corresponds to $\rho_0(\Gamma_1)*_{\langle\rho_0(\gamma)\rangle} O_s \rho_0(\Gamma_2)
O_s^{-1}$.}

A similar statement holds for the right component of $S\setminus\gamma$ if we choose a base point on the right component. {This case corresponds to $O_s^{-1}\rho_0(\Gamma_1)O_s*_{\langle\rho_0(\gamma)\rangle}  \rho_0(\Gamma_2)$.}

(2) If $\gamma$ is non-separating, the limit structure is the one obtained by attaching cylinders along each copy of $\gamma$.
\end{theorem}
\begin{proof}If $s\ra \infty$, then $C$, the boundary arc on
$\partial \Omega$ bounded by $\gamma_0,\gamma_+$, moves toward to
$\overline{\gamma_0, \gamma_-}\cup \overline{\gamma_0,\gamma_+}$.
Similarly if $s\ra -\infty$, $C$ moves toward to
$\overline{\gamma_-,\gamma_+}$.   If $\gamma$ is separating
one can see that as $s\ra\infty$, the Hausdorff limit domain $\Omega_\infty$ with a base point in
$K^L$ is equal to the one describe in Lemma \ref{1} and the Hausdorff limit domain with a base point in $K^R$ is $K^R$ itself as in Lemma \ref{2}.  For non-separating curve $\gamma$,
the Haudorff limit domain is obtained from $K^L$ or $K^R$ by attaching triangles or by cutting
out the domains along the edges of $K^L$ or $K^R$ respectively as in Lemma \ref{nonseparating}. When $s\ra-\infty$, the Hausdorff limit domain can be described similarly.   See Figure 1 and 2.

We also claim that $\Omega_\infty$ is contained in an affine chart. The only possibility that $\Omega_\infty$ is not
contained in an affine chart is that   vertices of  two triangles coincide, i.e. there is a triangle $\alpha_i(\triangle)$ whose
vertex $\alpha_i(\gamma_0)$ is at infinity. But vertex $\alpha_i(\gamma_0)$ is determined by two tangent lines at the end { points $\alpha_i(\gamma_\pm)$} of $\alpha_i(\tilde\gamma)$. This means that two tangent lines are parallel. One can change the affine chart so that one of end
points $\alpha_i(\gamma_\pm)$ can be moved a bit so that the tangent lines are no longer parallel. Note here that there are no sequence of lifts of $\gamma$ converging to $\alpha_i(\tilde\gamma)$. Hence by moving one point, we do not create  new parallel lines at the end points of some lift of $\gamma$. Then by using the convexity of the domain and
two dimensional geometry, and the fact that the lifts of $\gamma$ are disjoint, no two tangent lines at the end points
of a lift of $\gamma$ are parallel. { Otherwise, to get two parallel lines at the end points of a lift $\tilde\gamma'$, $\tilde\gamma'$ and $\alpha_i(\tilde\gamma)$ would intersect.} This way we can find an affine chart for all of $\Omega_s$ such that the limit domain $\Omega_\infty$ is still contained in this affine chart.

In terms of representation, if the
original holonomy representation is
$\rho_0(\pi_1(S))=\rho_0(\Gamma_1)*_{\langle \rho_0(\gamma)\rangle}\rho_0( \Gamma_2)$, then bulging
deformation gives rise to $\rho_0(\Gamma_1)*_{\langle\rho_0(\gamma)\rangle} O_s \rho_0(\Gamma_2)
O_s^{-1}$.  If we write the matrix in the basis of $\gamma_+,\gamma_0,\gamma_-$, and $$\alpha=\left (\begin{matrix}
                  a & b & c \\
                  d & e & f \\
                  g & h & i \end{matrix}\right) \in \rho_0(\Gamma_2) \subset
                  SL(3,\mathbb R),$$ then
$$O_s \alpha O_s^{-1}=\left (\begin{matrix}
             a & e^{-s}b &  c \\
             e^s d & e  & e^s f\\
             g &  e^{-s}h & i \end{matrix}\right).$$
Since $O_s \alpha O_s^{-1}$ acts as an element in $PGL(3,\mathbb R)$, when
$s=\infty$ either it acts as $$\left (\begin{matrix}
          a & 0 & c \\
          0 & e & 0 \\
          g & 0 & i \end{matrix}\right)$$ if $d=f=0$ (in this case it stabilizes the axis of $\gamma$, and since
          $\rho(\pi_1(S))$ is discrete, it should be the power of $\gamma$), or as
          $$ \left (\begin{matrix}
          0 & 0 & 0 \\
          d & 0 & f \\
          0 & 0 & 0 \end{matrix} \right).$$ In the first case it acts
          along the axis of $\gamma$, and in the second case  it
          maps everything to $\gamma_0$. For $s=-\infty$, it acts as
$$\left (\begin{matrix}
          0 & b & 0 \\
          0 & 0 & 0 \\
          0 & h &0 \end{matrix}\right),$$ hence it maps everything
          to the axis of $\gamma$.

This shows that $\Gamma_2$ part disappear in the limit, and the
projective structure is supported only on $S'$.

For   $\gamma$ nonseparating,
$$\rho(\pi_1(S))=\Gamma=\Gamma_1*_{\langle\gamma\rangle}=\langle \Gamma_1,
\gamma_2\rangle\,,$$
 where  $\gamma_2\in\Gamma\setminus\Gamma_1$ conjugates two subgroups of
$\Gamma_1$ that are isomorphic to $\langle\gamma\rangle$. The bulging
deformation gives rise to
$$\langle \Gamma_1,
O_s\gamma_2\rangle.$$
The similar analysis shows that either the
convex structure converges to the one with two half open cylinders
attached to the two copies of $\gamma$, or to the one with two
cylinders corresponding to $\gamma$.

For the limit projective structure on the right side of $\gamma$, one can do a similar
analysis with $$O_s^{-1}\rho_0(\Gamma_1) *_{\langle \rho_0(\gamma) \rangle} O_s\rho_0(\Gamma_2) O_s^{-1} O_s=O_s^{-1}\rho_0(\Gamma_1) O_s *_{\langle \rho_0(\gamma) \rangle}\rho_0( \Gamma_2).$$ \vskip .1 in


Such a Hausdorff limit can be seen as follows geometrically. {We explain it first when  $\gamma$ is separating.}
Let $\alpha\neq\gamma$ be a fixed simple geodesic loop in the left component $S'$ of $S\setminus \gamma$, $\beta\neq\gamma$ a fixed simple geodesic loop on the right component of $S\setminus\gamma$.
We show that the distance between $\alpha$ and $\beta$ goes to $\infty$ as $s\ra\infty$.
Let
$$\rho_0(\alpha)=\left (\begin{matrix}
                  a & b & c \\
                  d & e & f \\
                  g & h & i \end{matrix}\right),\ \rho_0(\beta)=\left (\begin{matrix}
                  b_1 & b_2 & b_3\\
                  b_4 &  b_5 & b_6 \\
                  b_7 & b_8 & b_9 \end{matrix}\right)$$ written  in the basis of $\gamma_+,\gamma_0,\gamma_-$.  Then

$$\rho_s(\alpha\beta)=\left(\begin{matrix}
                  a & b & c \\
                  d & e & f \\
                  g & h & i \end{matrix}\right) \cdot    \left (\begin{matrix}
                  b_1 & e^{-s}b_2 & b_3 \\
                 e^sb_4 &  b_5 & e^sb_6 \\
                  b_7 & e^{-s}b_8 & b_9 \end{matrix}\right)$$

Hence $\text{ Trace}(\rho_s(\alpha\beta))=e^s(bb_4+hb_6)+\cdots$.
If $bb_4+hb_6=0$ for any $\alpha$ in the left component, $b_4=b_6=0$.  Then
$\rho_0(\beta)$ will fix the invariant geodesic of $\gamma$. This implies that $\beta$ must be a power of $\gamma$, which is a contradiction. If $bb_4+hb_6=0$ for any $\beta$ in the right component, then
$b=h=0$. In this case, $\rho(\alpha)$ will fix  $\gamma_0$, hence its invariant geodesic
must coincide with that of $\gamma$, again a contradiction. Hence there must be a pair of geodesics
$\alpha$ and $\beta$ such that the trace of $\rho_s(\alpha\beta)$ is in the order of $e^s(bb_4+hb_6)$.
{ Note here that if the trace of a hyperbolic isometry is equal to $\lambda_1+\lambda_2+\lambda_3$, then the Hilbert length of the closed geodesic defined by that isometry
is equal to $\frac{1}{2}\log \frac{\lambda_1}{\lambda_3}>\frac{1}{2}\log (\lambda_1)$. Hence if the trace tends to $\infty$, so does the Hilbert length.}

Let $\alpha_s$ and $\beta_s$ be geodesic representatives of $\rho_s(\alpha)$ and $\rho_s(\beta)$.
Choose a geodesic segment $l_s$ connecting $\alpha_s$ and $\beta_s$ which realizes the distance between them. Then by considering the arc $\alpha_s\cup l_s\cup \beta_s$, we have
$$\ell(\rho_s(\alpha\beta))\leq \ell(\rho_s(\alpha))+ 2\ell(l_s)+ \ell(\rho_s(\beta)).$$
Here $\ell$ denotes the Hilbert length.
Since $\ell(\rho_s(\alpha\beta))\ra\infty$ and $\ell(\rho_s(\alpha))=\ell(\rho_0(\alpha)), \ell(\rho_s(\beta))=\ell(\rho_0(\beta))$, we have $\ell(l_s)\ra\infty$.

{For non-separating $\gamma$, the extra generator $\gamma_2$ for $\rho(\pi_1(S))=\langle \Gamma_1,\gamma_2\rangle$ can be chosen as the image of a closed curve $\alpha$ intersecting
$\gamma$ at one point.  When we write $\gamma_2=\begin{pmatrix}
 a & b & c\\
  d & e & f \\
  g & h & i \end{pmatrix}$ in the basis $\gamma_+,\gamma_0,\gamma_-$, by choosing $\gamma_2$ properly, we can assume that $e\neq 0$ and {$a+i\neq 0$}. Under bulging deformation, it is changed to
$$O_s\gamma_2 =\begin{pmatrix}
                  e^{-\frac{s}{3}}a & e^{-\frac{s}{3}}b & e^{-\frac{s}{3}}c \\
                 e^{\frac{2s}{3}} d & e^{\frac{2s}{3}} e & e^{\frac{2s}{3}} f \\
                  e^{-\frac{s}{3}}g & e^{-\frac{s}{3}} h & e^{-\frac{s}{3}}i \end{pmatrix}.$$}
The trace is $e^{-\frac{s}{3}}(a+i)+e^{\frac{2s}{3}}e=e^{-\frac{s}{3}}(a+e+i)+(e^{\frac{2s}{3}}-e^{-\frac{s}{3}})e$.
Then the trace of $O_s\gamma_2$ tends to infinity as $s\ra\pm\infty$. Hence the length of the transversal curve $\alpha$ to $\gamma$ tends to infinity.

{Since above argument is true for any curve on the left component and any curve on the right component for $\gamma$ separating, and for transversal curve to $\gamma$ when $\gamma$ is non-separating, there exists a cylinder containing the geodesic $\gamma$ connecting two components.}
More precisely
\begin{lemma}\label{l}As $s\ra\infty$, the length of the cylinder between left and right component of $S\setminus\gamma$ grows in the order of $s$.
\end{lemma}
\begin{proof}Above argument shows that the distance between left and right component grows in the order of $\ell(\rho_s(\alpha\beta))$. The trace $\lambda_1+\lambda_2+\lambda_3, (\lambda_1>\lambda_2>\lambda_3, \lambda_1\lambda_2\lambda_3=1)$ grows in the order of $e^s$ along the bulging deformation, and the Hilbert length $\frac{1}{2}\log \frac{\lambda_1}{\lambda_3}=\frac{1}{2}\log (\lambda_1^2\lambda_2)$
grows in the order of $\log \lambda_1$, hence the length of the cylinder grow in the order of
$s$.
\end{proof}
Take a fundamental domain $F_s$ of $\rho_s(\pi_1(S))$ in $\Omega_s$ which is divided by a segment of $\tilde\gamma$ into left side $F_s^L$ and right side $F_s^R$. Then the diameter
of $F_s$ tends to $\infty$ by the above argument.  Since $\rho_s(\pi_1(S))=\rho_0(\Gamma_1)*_{\rho_0(\gamma)}O_s\rho_0(\Gamma_2)O_s^{-1}$, the action of $\rho_0(\Gamma_1)$ is always the same on $K^L$. By choosing the fundamental domain $F_s$
inside $K^L\cup O_s(K^R)\subset \Omega_s$ by  Equation (\ref{sepL}), we see that $F_s^L=F_0^L$ and $F_s^R=O_s(F_0^R)$.  By considering $\cup_{n\in\mathbb Z} \rho_s(\gamma)^n(F_s)$, we can see that the right side $\cup_{n\in\mathbb Z}\rho_s(\gamma)^n(F_s^R)$ tends
to $\triangle$. Considering this process for each lift of $\gamma$, we see that
the Hausdorff limit with a base point on the left side of $S\setminus\gamma$ converges to the one described in Lemma \ref{1}. This finishes the proof for the claim in the statement.
\end{proof}

{For any $s$, since $\gamma$ has the same length and $\gamma$ is the core curve of the cylinder,}
 the injectivity radius  $r_s$ along $l_s$ is bounded below independent of $s$ for large $s$ by the length of $\gamma$. { Consider the band $l_s\times L_s$, where
the length of $L_s$ is $r_s$, and $\ell(l_s)$ goes to infinity. Since the Hausdorff measure is induced by the Hilbert metric, the volume of this band tends to infinity.} This shows that
\begin{theorem}\label{Area}Let $\rho_s:\pi_1(S)\ra SL(3,\mathbb R)$ be a family of bulging deformation along a simple closed curve $\gamma$. Then the Hilbert area of the corresponding real projective surface has
$$\text{Area}(\Omega_s/\rho_s(\pi_1(S)))\ra\infty.$$
\end{theorem}

As  in $\bc P^1$-structure case, this construction generalizes
to any measured lamination with transverse measure in $\mathfrak a$, the Lie algebra of $A$ in the Iwasawa decomposition $KAN$ of $SL(3,\mathbb R)$.
This construction combines the earthquake deformation by Thurston
and bulging deformation. Since the earthquake deformation gives
$6g-6$ dimensional deformation, and bulging also gives the same dimensional deformation, near
the Teichm\"uller space, this construction gives $6g-6$ dimensional
deformation of the Teichm\"uller space, {hence $12g-12$ dimensional subset   near Teichm\"uller space in the space of real convex structures.}  Missing $4g-4$ dimensional
deformation comes from the fact that the projective structures on
a pair of pants are parametrized by 2-dimensional parameters for
each boundary together with 2 more parameters inside. Since there
are $2g-2$ pairs of pants, there are $2(2g-2)$ parameters to
determine the projective structures inside pairs of pants.

It is known that the set of convex projective structures on a closed surface is homeomorphic to the set of pairs $(\Sigma, U)$ where $\Sigma$ is a hyperbolic surface and $U$ is a holomorphic cubic differential over $\Sigma$ \cite{La, Lo2}.
 In
\cite{Lo}, it is shown that if $(S_i,U_i)$ is a sequence such that
$S_i$ tends to a nodal curve pinched along $\gamma$ at the boundary
of the Deligne-Mumford compactification of the moduli space of
Riemann surface, and the cubic differential $U_i$ tends to a regular
cubic differential $U$ over the nodal curve, then the vertical twist
parameter along $\gamma$ tends to $\pm\infty$ if the residue $R$ of
$U$ along $\gamma$ is nonzero and $Re R\neq 0$.

Using this, one can prove the following.
\begin{proposition}\label{infinite}If a sequence of strictly convex projective
structures converge to a non-strictly convex projective structure,
then {the area of the limit structure is infinite}. A concrete example can be obtained
by bulging deformation along a simple curve $\gamma$. In this
example, the projective surface converge to a projective surface  of
infinite area with a  boundary corresponding to $\gamma$, and the
cubic differential converges to a regular cubic differential over
this limit surface.
\end{proposition}
\begin{proof}
If the limit convex projective structure has a finite
covolume,  the limit domain should be strictly convex \cite{Marquis}.
In bulging
deformation, if we let $s\ra \pm\infty$, the domain will converge to
the one with  line segments on the boundary by Proposition \ref{cylinder}. In either case, the holonomy representing
the puncture is hyperbolic and the surface is of infinite volume, {which is stated in  Theorem \ref{Area} with a different method.}

In \cite{Lo} (Theorem 3), it is shown that the holonomy type and
vertical twist parameters vary continuously as long as the regular
cubic differentials at the limit has non-vanishing residue. But if
the holonomy of the puncture  is hyperbolic, then the residue is
nonzero. Indeed, for $s\rightarrow \infty$, the limit structure has a triangular end, i.e. quotient of a triangle by a hyperbolic isometry, whose three vertices are
three fixed points of the hyperbolic isometry, and the cubic differential has a third order pole with $Re(R)<0$. For $s\rightarrow -\infty$, the limit domain contains
the axis of the hyperbolic isometry on the boundary and the cubic differential has a third order pole with $Re(R)>0$. In either case, the cubic differential has a third order pole with $Re(R)\neq 0$.
See  \cite{Nie2} for details.
 Hence the
cubic differentials converge to a regular cubic differential over
this projective surface.
\end{proof}


\section{Entropy of Bowen-Margulis measure, topological entropy}
It is a well-known fact that the topological entropy is equal to the entropy of the Bowen-Margulis measure $\mu_{BM}$ for
strictly convex real projective structure, see \cite{Crampon}.

Consider the deformation $\rho_{s}=\rho(\Gamma_1)*_{\langle \gamma \rangle} ( O_s)\rho(\Gamma_2)( O_s)^{-1}$ or $\langle \rho(\Gamma_1),
O_s\gamma_2\rangle$ as $s\ra\infty$ where $\pi_1(S)=\Gamma_1*_{\langle \gamma \rangle} \Gamma_2$ or HNN extension depending on whether $\gamma$ is separaing or not. For any $\alpha\in\Gamma_1$, through the deformation, the Hilbert length of a closed geodesic corresponding to $\alpha$ is fixed by $\frac{1}{2}\log \frac{\lambda_1}{\lambda_3}$ where $\lambda_1>\lambda_2>\lambda_3$ are
the eigenvalues of $\rho_{0}(\alpha)$. Also for $( O_s)\rho(\Gamma_2)( O_s)^{-1}$, since the representation is just a conjugation, the closed geodesic lengths remain the same.

But for any $\alpha*\beta$ such that $\alpha\in\Gamma_1,\beta\in\Gamma_2$,  the length goes to infinity as $s\ra\infty$ by Lemma \ref{l}. Hence it is {believed} that the topological entropy converges to a positive
number greater than the topological entropy of $\rho(\Gamma_1)$ and $\rho(\Gamma_2)$.

To prove the convergence of the topological entropy {as $s\ra\infty$}, we describe the combinatorial method of a closed curve.
Take a pants decomposition {$P=\{\gamma_1,\cdots,\gamma_{3g-3}\}$} including $\gamma=\gamma_1$. For each pair of pants in $S\setminus P$, take an ideal triangulation into two ideal triangles.  Each closed curve $\alpha\notin P$
intersects this ideal triangulation $\cal T$  transversely. Note that the edge of the triangulation accumulates on two distinct curves in $P$.
Along the orientation of $\alpha$, one can define a map $$suc:\alpha\cap\cal T\ra \alpha\cap\cal T$$ which maps $p$ to the next intersection point.

 Call a point $p$ a crossing point if there is no $\beta\in P$ such that all the edges of $\cal T$ containing $p$, $suc^{-1}(p)$, $suc(p)$ accumulate on $\beta$. Call a segment $[suc^{-1}(p), suc(p)]\subset\alpha$ a {\bf crossing segment} of $\alpha$ at $p$. The edge containing $p$ is the one running from one curve in $P$ to another curve in $P$ in the same pair of pants. Then combinatorially there are at most 6 different un-oriented crossing segments in a pair of pants, { 3 in the front, 3 in the back}, which gives at most $12(2g-2)$ different oriented crossing segments on a closed surface of genus $g$.

 These crossing segments are connected by two types of segments of $\alpha$, one is pants-changing and the other looping segment.
 The {\bf pants-changing segment} is the one changing the pants, and the {\bf looping segment} is the one staying in one pants but looping around the boundary of the pants.

To measure how many times the looping and pants-changing segment winds around the curve $\gamma_i\in P$, we introduce the following quantities.  Let $P_1$ and $P_2$ be two pairs of pants glued along $\gamma_i$. Let $l_1\subset P_1$ be a leaf of the ideal triangulation which does not accumulate on $\gamma_i$, $l_2\subset P_2$ which does not accumulate on $\gamma_i$. Choose the length minimizing arc $d_i$ connecting $l_1$ and $l_2$ which intersects $\gamma_i$ transversely.

If $\beta$ is a pants-changing segment, let $\#(\beta)$ denote either $0$ if $d_i$ is not transverse to $\beta$ or $|\beta\cap d_i|=0$, or $|\beta\cap d_i|-1$ when $|\beta\cap d_i|$ is
the positive number of intersection points. Obviously $\#(\beta)$ measures the number of times that $\beta$ winds around
$\gamma_i$.  Similarly if $\beta$ is a looping segment winding around $\gamma_i\in P$, let $\#(\beta)$, { with the same notation,} denote the number of self intersections of $\beta$. It also measures the number of times that $\beta$ winds around $\gamma_i$.

Suppose $m$-crossing segments are given. Then one can form a closed loop by connecting these crossing segments by pants-changing or looping segments $\beta$. Suppose $\#(\beta)$ are given also. Then it is not difficult to show that there are at most $18^m$ different closed loops with these data,
{ where $m$ is the number of crossing segments}. See \cite{Zhang} (Proposition 3.22).

Now we want to estimate the number of closed loops whose length is less than $T$ using these combinatorial data.

{\bf Suppose $\eta$ is a closed geodesic which is not homotopic to any curve in $S\setminus P$}. Then using the orientation of $\eta$, there are natural orderings on crossing points, crossing segments and pants-changing, looping segments.
{When we do bulging deformation along a simple closed curve $\gamma=\gamma_1$}, the Hilbert geometry on each pair of pants remains bounded below except  $P_1$ which is bounded by $\gamma=\gamma_1$ in the following sense.
\begin{lemma}Along bulging deformation on $\gamma=\gamma_1$, there exists a lower bound for the length of crossing segments in every pairs of  pants.
\end{lemma}
\begin{proof}A crossing segment in the lift of the ideal triangulation on a pair of pants is the one joining two edges of ideal triangles as in Figure 3 where $\alpha,\beta,\delta$ are boundary components of the pair of pants.

\begin{figure}
\scalebox{1.1} 
{
\begin{pspicture}(1,-0.238125)(17.881874,2.238125)
\psellipse[linewidth=0.04,dimen=outer](7.08,0.6296875)(1.66,1.47)
\psline[linewidth=0.04cm](7.44,2.07)(8.7, 1.0)
\psline[linewidth=0.04cm](7.44,2.07)(6.3, 1.9)
\rput(5.9,1.9){$\tilde\beta$}
\psline[linewidth=0.04cm](8.7, 1.0)(8.5,-0.05)
\psline[linewidth=0.04cm](8.5,-0.05)(7.3,-0.83)
\psline[linewidth=0.04cm](8.5,-0.05)(8.7,0.6)
\rput(8.95,0.4){$\tilde\alpha$}
\psline[linewidth=0.04cm](7.3,-0.83)(5.6,-0.05)
\psline[linewidth=0.04cm](5.6,-0.05)(6.8,-0.8)
\rput(6,-0.9){$\tilde \delta$}
\psline[linewidth=0.02cm](5.7,1)(8.02,1)
\usefont{T1}{ptm}{m}{n}
\psline[linewidth=0.04cm](5.6,-0.05)(5.8,1.5)
\psline[linewidth=0.04cm](5.8,1.5)(7.44,2.07)
\psline[linewidth=0.04cm](7.44,2.07)(8.5,-0.05)
\psline[linewidth=0.04cm](8.5,-0.05)(5.6,-0.05)
\psline[linewidth=0.04cm](5.6,-0.05)(7.44,2.07)
\end{pspicture}
}
\caption{A crossing segment intersecting the lift of  ideal triangulation. Three triangles around the center one are identified under the action of $\alpha,\beta,\delta$}
\end{figure}
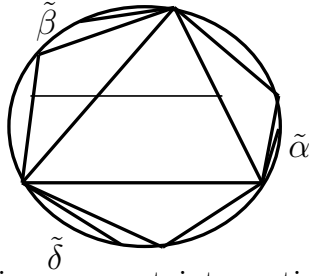

Suppose a lift of the crossing segment is on the left side of $\tilde\gamma$, hence contained in $K^L$.  By Equations (\ref{sepL}) and (\ref{nonsep}),  $K^L\subset\Omega_s\subset \Omega_\infty^L$.
By Section \ref{Hilbertmetric} the length of crossing segment will be  bounded  below by the length determined in $\Omega_\infty^L$.
Note that since  $\rho_s= \rho(\Gamma_1)*_{\langle \gamma \rangle} O_s \rho(\Gamma_2) O_s^{-1}$ (or $\langle \rho(\Gamma_1), O_s\gamma_2\rangle$),  the action of $\Gamma_1$ on $K^L$ is the same for all $s$.  Since the  action of $\Gamma_1$ on $K^L$ has a compact fundamental domain $F$ {where the lift of the crossing segment lies in}, there is a lower bound for the length of the crossing segments
{measured in $\Omega_\infty^L$.}
When the lift of the crossing segment is in $K^R$, {we apply the same logic to $O_s^{-1}\rho(\Gamma_1)O_s*_{\langle \gamma \rangle}\rho(\Gamma_2)$.}  In  either case, the length of the lift of the crossing segment is bounded below.
\end{proof}
 Since the geometry of $P_1$ is the one $P_1'$ with a long cylinder $\cal C$ attached to $\gamma_1=\gamma$ by Proposition \ref{cylinder}, we may use  pairs of pants $P_1',P_2,\cdots,P_{2g-2}$ and the cylinder $\cal C$ of height $s$ (Lemma \ref{l}) attached to $\gamma_1$ to estimate the length of a closed geodesic. Note that the length of crossing segments on these new pants are uniformly bounded below and the boundary components of these pants have uniformly bounded lengths throughout the bulging deformation along $\gamma_1$. If $\eta$ intersects $\gamma_1$, the length of all crossing
segments are uniformly bounded below by $Cr$ independent of $\eta$ on each pair of pants $P_1',P_2,\cdots,P_{2g-2}$, and just add $s$ to the pants-changing segment which crosses $\gamma_1$. Also it is clear that the length
of each looping or pants-changing segment $\beta$,
$$\ell_s(\beta)\geq \#(\beta) L,$$ where $L$ is $\min\{\ell_s(\gamma_i)=\ell_0(\gamma_i)|\gamma_i\in P\}$, a fixed number for all $s$.
Hence for any $\eta$ intersecting $\gamma=\gamma_1$,
\begin{eqnarray}\label{length}
\ell_s(\eta)\geq m Cr + \sum_\beta \#(\beta) L +s i(\eta,\gamma)
 \end{eqnarray} where the sum is over all pants-changing or looping segments, and $m$ is the number of crossing segments.
We can prove the following theorem for a bulging deformation performed on  every curve of pants decomposition. {One should understand that performing bulging deformation on $P$ is doing bulging deformation on each curve in $P$ one by one.}
\begin{theorem}\label{entropyzero}Let $P$ be a pants decomposition of $S$. Then the topological entropy of $\rho_s$ converges to $$\max_{P_i} \lim_{T\ra\infty} \frac{\log \#\{[\eta]|\ell_s(\eta)=\ell_0(\eta)<T, \eta\in  P_i \}}{T}   $$ as $s\rightarrow \infty$ where $\rho_s$  is a bulging deformation on $P$ and $P_i$ are components of $S\setminus P$.
\end{theorem}
\begin{proof}Since we do bulging along each curve in $P$, the resulting one has cylinder of width at least $\frac{s}{3}$ for each curve in $P$. By choosing $\gamma_i$ in the middle of cylinder, we may assume that each pair of pants has cylinder neighborhood of each boundary component of width $\frac{s}{6}$. By modifying ideal triangulation if necessary, we may assume that the part of leaves of ideal triangulation which wrap around the curves in $P$ occur only in the cylinder regions of width $\frac{s}{12}$. { Note here that neither the pants decomposition curves, nor ideal triangulation curves are geodesics in general.} Hence $Cr_s$, the minimum length of crossing segments, tends to $\infty$ as $s\ra\infty$.
Then for any $\eta$ not in $S\setminus P$, the equation (\ref{length}) becomes
\begin{eqnarray}\label{length2}
\ell_s(\eta)\geq m Cr_s + \sum_\beta \#(\beta) L
 \end{eqnarray} where $L$ is fixed independent of $s$ due to Proposition \ref{cylinder}. Then
 $$m\leq \frac{1}{Cr_s}\ell_s(\eta).$$ 

{\bf If $\eta$ is contained in $S\setminus P$, then when we represent $\eta$ as a geodesic, we cannot apply
above estimate since the length of the crossing segment might be bounded}.  Indeed this happens for all curves entirely contained in one of pairs of pants when we do bulging along $P$ since the representations
restricted to a pair of pants are just conjugations of the original ones.
{\bf Now we count the number of geodesics whose length is less than $T$ and not lying entirely in $S\setminus P$}.
$$\#\{[\eta]|  \ell_s(\eta)<T\}<\#\{[\eta]|m Cr_s + \sum_\beta \#(\beta) L <T\} $$
$$\leq \sum_{m=1}^{\lfloor\frac{T}{Cr_s}\rfloor} 18^m \#\{m Cr_s + \sum_\beta \#(\beta) L <T\}  $$
$$\leq \sum_{m=1}^{\lfloor\frac{T}{Cr_s}\rfloor} 18^m (24g-24)^m \#\{(\beta)|m Cr_s + \sum_\beta \#(\beta) L <T\} $$
$$ \leq \sum_{m=1}^{\lfloor\frac{T}{Cr_s}\rfloor} (432g-432)^{\lfloor\frac{T}{Cr_s}\rfloor} \#\{(\beta)|m Cr_s + \sum_\beta \#(\beta) L <T\} .$$
Here $\lfloor a \rfloor$ denotes the largest integer less than $a$.
The above inequalities follow since once a sequence of  $m$-crossing segments and  $\#(\beta)$ are given, there exist at most $18^m$ different curves, and since there are at most $(24g-24)$ different oriented crossing segments, there are at most $(24g-24)^m$
possible sequence for each positive $m$. Now if $m$ crossing segments are given, we need to connect them by looping or pants changing segments satisfying $k=\sum_\beta \#(\beta)  <\lfloor\frac{1}{L}(T- m Cr_s)\rfloor$. Hence we need to calculate the number $f_s(m,T)$ of ways to partition all integers $k<\lfloor\frac{1}{L}(T- m Cr_s)\rfloor$ into $m$ non-negative integers. Then
the above inequality is
$$\leq (432g-432)^{\lfloor\frac{T}{Cr_s}\rfloor} \sum_{m=1}^{\lfloor\frac{T}{Cr_s}\rfloor} f_s(m,T).$$
Hence the topological entropy is less than
$$\lim_{T\ra\infty} \frac{1}{T}[\frac{T}{Cr_s}\log(432g-432)+\log \sum_{m=1}^{\lfloor\frac{T}{Cr_s}\rfloor} f_s(m,T)].$$
But as $\frac{1}{T}\frac{T}{Cr_s}\ra 0$ as $s\ra\infty$, it suffices to show that
$$\lim_{s\ra\infty}\lim_{T\ra\infty} \frac{1}{T}[\log \sum_{m=1}^{\lfloor\frac{T}{Cr_s}\rfloor} f_s(m,T)]=0.$$
Now we need to find out $f_s(m,T)$. The number of ways to partition $k$ into $m$ non-negative integers is equal to the number of non-negative integer solutions to $x_1+\cdots+x_m=k<\frac{1}{L}(T- m Cr_s)$. Hence
$$f_s(m,T)\leq\sum_{k=0}^{\lfloor\frac{1}{L}(T- m Cr_s)\rfloor} {m+k-1 \choose k}$$ where ${m+k-1 \choose k}$ denotes the number of ways to  choose $k$ elements out of $m+k-1$ elements. Suppose that $f_s(M_s,T)$ is maximum among $f_s(m,T), 1\leq m\leq \lfloor\frac{T}{Cr_s}\rfloor$ and $ {M_s+q_s-1 \choose q_s}$ is maximum among $ {M_s+k-1 \choose k}, 0\leq k\leq \lfloor\frac{T-M_s Cr_s}{L}\rfloor$. Then
$$\lim_{s\ra\infty}\lim_{T\ra\infty} \frac{1}{T}[\log \sum_{m=1}^{\lfloor\frac{T}{Cr_s}\rfloor} f_s(m,T)]\leq $$
$$\lim_{s\ra\infty}\lim_{T\ra\infty} \frac{1}{T}\log(\frac{T}{Cr_s}\frac{T-M_s Cr_s}{L}  {M_s+q_s-1 \choose q_s})$$$$\leq
\lim_{s\ra\infty}\lim_{T\ra\infty} \frac{1}{T}\log {M_s+q_s \choose q_s}.$$ One can show that this last thing converges to zero as follows.
The Stirling's Formula tells that
$$n!\sim (\frac{n}{e})^n\sqrt{2\pi n}.$$ Hence
$${M_s+q_s \choose q_s}=\frac{(M_s+q_s)!}{M_s!q_s!} \sim \frac{1}{\sqrt 2\pi}\sqrt{\frac{M_s+q_s}{M_s q_s}}\frac{(M_s+q_s)^{M_s+q_s}}{M_s^{M_s}q_s^{q_s}}. $$
Now it is easy to see that
$$\lim_{T\ra\infty}\frac{1}{T}\log \sqrt{\frac{M_s+q_s}{M_s q_s}}\ra 0$$ using $ q_s\leq \frac{T}{L},\ M_s\leq \frac{T}{Cr_s}$.
The second term
$$\lim_{T\ra\infty}\frac{1}{T}\log \frac{(M_s+q_s)^{M_s+q_s}}{M_s^{M_s}q_s^{q_s}}=\lim_{T\ra\infty}\frac{1}{T}\log
(1+\frac{q_s}{M_s})^{M_s}(1+\frac{M_s}{q_s})^{q_s}.$$
If $\frac{M_s}{q_s}\ra\alpha\neq 0$, then since
$$\frac{M_s}{T}\leq \frac{1}{Cr_s}$$
the above quantity is less than
$$\frac{1}{Cr_s}\log(1+\alpha^{-1})+\frac{1}{\alpha Cr_s}\log(1+\alpha)),$$ which converges to zero as $s\ra\infty$.
If $\frac{M_s}{q_s}\ra 0$, then it is easy to show that the limit is zero.

If $\frac{M_s}{q_s}\ra\infty$, then
$$\lim_{T\ra\infty} \frac{1}{T}\log(1+\frac{q_s}{M_s})^{M_s}(1+\frac{M_s}{q_s})^{q_s}$$$$=\lim_{T\ra\infty} [\frac{q_s}{T}\frac{M_s}{q_s}\log(1+\frac{q_s}{M_s})+\frac{M_s}{T}\frac{q_s}{M_s}\log(1+\frac{M_s}{q_s})]\leq$$
$$\lim_{T\ra\infty}\frac{q_s}{T}+\lim_{T\ra\infty}\frac{1}{Cr_s} \frac{q_s}{M_s}\log(1+\frac{M_s}{q_s})=0 $$ using
$\lim_{x\ra 0}\frac{\log(1+x)}{x}=1$.

Let $S\setminus P=P_1\cup\cdots \cup P_{2g-2}$.
Then the topological entropy of $\rho_s$ can be estimated by
$$\lim_{T\ra\infty} \frac{\log[ \#\{[\eta]|\ell_s(\eta)=\ell_0(\eta)<T, \eta\in S\setminus P\}+\#\{[\eta]|\ell_s(\eta)<T, \eta\notin S\setminus P\}]}{T}  $$
$$\leq \lim_{T\ra\infty} \frac{\log \#\{[\eta]|\ell_s(\eta)=\ell_0(\eta)<T, \eta\in S\setminus P\}}{T}$$
$$\leq \max_{P_i} \lim_{T\ra\infty} \frac{\log \#\{[\eta]|\ell_s(\eta)=\ell_0(\eta)<T, \eta\in  P_i \}}{T}.    $$    The first inequality follows from the fact 
$$\lim_{T\ra\infty} \frac{\log\#\{[\eta]|\ell_s(\eta)<T, \eta\notin S\setminus P\}}{T}=0.$$
This shows that the topological entropy of $\rho_s$ converges to
$$\max_{P_i} \lim_{T\ra\infty} \frac{\log \#\{[\eta]|\ell_s(\eta)=\ell_0(\eta)<T, \eta\in  P_i \}}{T}.    $$
\end{proof}
\begin{corollary}For each $\alpha\in[0,1]$, there exists a divergent sequence in the space of convex real projective structures whose topological entropy converges to $\alpha$.
\end{corollary}
\begin{proof}Fix  $P$ a pants decomposition. On a pair of pants, it is not difficult to show that there exists a projective structure with geodesic boundary that the topological entropy is less than any designated number near zero. Indeed the estimate (\ref{length2}) and the subsequent arguments hold if we can gurantee that $Cr_s\ra\infty$ as $s\ra\infty$. Such a sequence can be produced by letting the internal Goldman parameters diverge. For our case since bulging deformation changes only bulging parameters, the lengths of curves contained in $S\setminus P$ remain the same throughout the bulging deformation. Hence to increase the lengths of curves in $S\setminus P$, we have to let the internal parameters diverge.
Indeed it is shown in \cite{Zhang} that the shortest length goes to infinity as internal parameters diverge.

 Produce projective structures $X_i$ by gluing pairs of pants with small topological entropy  so that the maximum of topological entropy on each pair of  pants $S\setminus P$ is  $\alpha_i\ra 0$. Then do the bulging along $P$ to obtain a sequence
$X_i^s$. By Theorem \ref{entropyzero}, the topological entropy of $X_i^s$ goes to $\alpha_i$ as $s\ra\infty$. Let $X_0$ be a hyperbolic surface. Connect $X_0$ to $X_i$ by a continuous path. This path together with a bulging deformation $X_i^s$ from $X_i$ gives a continuous family of deformations starting from a hyperbolic surface. Since the topological entropy is a continuous function with respect to a continuous parameter, by taking a diagonal sequence, the topological entropy varies from 1 to $\alpha_i$, and hence there exists some $i$ and $s$ such that the topological entropy of $X_i^s$ is $\alpha$.
\end{proof}

When we do bulging deformation along a separating simple closed curve, it is believable that the topological entropy converges to a positive
number greater than the topological entropy of each side. Yet  it is not obvious to show that \\
{\bf Conjecture}:\ The topological entropy of $\rho_{t_0,s}$ decreases and converges to  the maximum of the topological entropy
of $\rho_{t_0,0}|_{\Gamma_1}$ and the topological entropy of  $\rho_{t_0,0}|_{\Gamma_2}$.
\vskip .3 in
{\bf Acknowledgement}\ We thank Tengren Zhang for pointing out some mistakes in the first draft of the paper. The  second author also thanks CIRM (Centre International de Rencontres Math\'ematiques) for arranging the visit to the first author.

\vskip .2 in
\noindent     Patrick Foulon\\ Institut de Recherche Mathematique Avanc\'ee,
UMR 7501 du Centre National de la Recherche Scientifique, 7 Rue Ren\'e Descartes, 67084
Strasbourg Cedex, France\\
\texttt{foulon\char`\@math.u-strasbg.fr}\\
\vskip .1 in
\noindent     Inkang Kim\\
     School of Mathematics\\
     KIAS, Heogiro 85, Dongdaemen-gu\\
     Seoul, 130-722, Korea\\
     \texttt{inkang\char`\@ kias.re.kr}
\end{document}